\documentclass[11pt]{article}
\usepackage[top=0.80in,bottom=0.80in,left=0.95in,right=0.95in]{geometry}
\usepackage[T1]{fontenc}
\usepackage{lmodern,amsmath,amssymb,amsthm,booktabs,microtype,graphicx,array,float,needspace}
\usepackage[colorlinks=true,allcolors=blue]{hyperref}
\usepackage{xurl}
\newtheorem{theorem}{Theorem}[section]
\newtheorem{proposition}[theorem]{Proposition}
\newtheorem{lemma}[theorem]{Lemma}
\newcommand{\EE}{\mathbb E}
\newcommand{\argmin}{\mathop{\rm arg\,min}}
\newcommand{\early}{\mathrm{early}}
\newcommand{\late}{\mathrm{late}}
\title{When later capacity losses require more inventory:\\a three-period classification with two-point demand}
\author{%
{\large Qichao Wang}\\[2pt]
{\normalsize Hebei University of Technology}\\[5pt]
{\large Daoyu Dong}\\[2pt]
{\normalsize University of Electronic Science and Technology of China}}
\date{}
\hypersetup{pdftitle={When later capacity losses require more inventory: a three-period classification with two-point demand},pdfauthor={Qichao Wang and Daoyu Dong}}
\begin{document}
\raggedbottom
\maketitle
\begin{abstract}
An earlier production-capacity loss can raise optimized inventory costs while requiring less initial stock than an equally large later loss. For a three-period backlogged inventory model with two-point demand and baseline future capacities $(\eta,2\eta)$, we identify six explicit cost thresholds that characterize the optimal target sets for every demand probability and positive cost ratio. In both open regions that violate the temporal-ordering conjecture of Jak\v{s}i\v{c} et al., we prove that the earlier loss has strictly higher optimized cost. A separate convex-operator argument establishes a general one-coordinate sensitivity bound, independently of the conjecture, and the classification attains the bound in both temporal directions. Strict reversals persist under small continuous demand perturbations. In a fixed-loss information experiment, imposing the conjectured ordering eliminates the entire value of advance loss-location information in the reversal regions.
\end{abstract}


\section{Introduction}
Advance capacity information (ACI) enables a manufacturer to adjust current inventory to future replenishment constraints. It is natural to expect that a loss of capacity in the near future calls for more precautionary stock than the same loss further ahead: an earlier disruption leaves less time to prepare. Yet the cost of a disruption and the marginal benefit of stocking an additional unit are different quantities. The distinction matters for planning rules that use the timing of a known capacity shortfall to rank inventory targets: a rule enforcing the wrong order can prevent a planner from choosing the optimal target in both scenarios.

A formal expression of the temporal intuition appears in Jak\v{s}i\v{c} et al.\ \cite[Conjecture 1]{Journal}: an earlier equal capacity reduction should produce at least as high a current base-stock target as a later reduction. The discussion preceding the conjecture further predicts equality when holding cost $h$ exceeds backlog cost $b$.\footnote{The temporal comparison appears as Theorem 4(1) in the working paper \cite{WP} and as Conjecture 2.1 in the thesis \cite{Thesis}.} A related claim, Conjecture 2 in \cite{Journal}, bounds the increase in the current target caused by reducing a single future capacity: it should be nonnegative and no greater than the capacity removed. These are distinct questions about the direction and the magnitude of inventory adjustments.

Consider three periods with independent, identically distributed demand equal to 0 or 2 with probability $1/2$ each, holding cost 1, backlog cost $1/2$, and baseline future capacities $(1,2)$. Removing one unit of capacity in period 2 produces a unique first-period target of 0 and an expected cost of $9/4$. Moving the loss to period 3 instead produces a target of 1 and a lower expected cost of $33/16$. The earlier disruption is more expensive even though the later disruption calls for more initial stock. This example also has $h>b$, so it contradicts both temporal assertions.

We develop four results. First, for every positive cost ratio and nondegenerate demand probability in the three-period family with baseline capacities $(\eta,2\eta)$, we characterize both target sets, including all boundary minimizer intervals. We also prove that throughout the two strict reversal regions the earlier loss is more expensive, showing that higher optimized cost need not imply a higher target. Second, because the published derivation of the one-coordinate sensitivity bound invokes Conjecture 1, we give an independent proof and establish its sharpness in both temporal directions. Third, we show that the two explicit reversals survive small continuous perturbations of the demand law. Finally, we quantify the cost of imposing the conjectured ordering on decisions that could otherwise use the loss location, and identify when that restriction eliminates the value of observing the loss location.

The structure of optimal inventory policies under production constraints is established in the capacitated inventory literature, including Federgruen and Zipkin \cite{FZ} and Kapu\'sci\'nski and Tayur \cite{KT}. Whereas these results concern the form of an optimal policy, our comparison concerns the ordering of its target levels across equally severe capacity losses. The value of advance information has also been studied from the demand side, where future demand interacts with limited production capacity \cite{OW}, and from the supply side, where forecasts evolve or disruptions alter feasible replenishment \cite{Altug,Atasoy}. In models with regenerative capacity \cite{Senga}, unused production capacity can affect later feasibility. In contrast, our period-specific capacities expire, making the timing of an equal loss an explicit comparative-statics variable. We ask how the distribution of future replenishment opportunities changes the marginal value of current stock. The information analysis in Section~\ref{sec:information} complements this literature by distinguishing the value of observing a loss location from the value that can actually be realized by a policy class constrained by a proposed structural property.

Section~\ref{sec:model} introduces the model and exact examples. Section~\ref{sec:main} gives the classification and its economic mechanism, followed by the sensitivity and stability results in Section~\ref{sec:sensitivity}. Section~\ref{sec:information} quantifies the cost of imposing the conjectured order and the value of loss-location information; Section~\ref{sec:conclusion} concludes. Appendix~\ref{app:slopes} contains the Bellman slope calculations.

\section{Model and motivating counterexamples}\label{sec:model}
\subsection{Model formulation}
There are three periods, zero lead time, nonnegative real orders, full backlogging, no ordering cost, no terminal salvage and no discount. End-of-period loss is
\[
 \ell(x)=h x_++b(-x)_+,\qquad h,b>0,
\]
where $x_+=\max(x,0)$. Demands are independent and identically distributed, with $D_t=2\eta B_t$, where $B_t$ are Bernoulli random variables with $\Pr(B_t=1)=p\in(0,1)$ and $\eta>0$. An order in period $t$ is bounded by known capacity $q_t$. For $V_4=0$, define
\begin{equation}\label{eq:bellman}
 J_t(y)=\EE[\ell(y-D_t)+V_{t+1}(y-D_t)],\qquad
 V_t(x)=\min_{0\le z\le q_t}J_t(x+z).
\end{equation}
The target set is $S_t=\argmin_{y\in\mathbb R}J_t(y)$, and we select its smallest element $s_t=\min S_t$. Targets are unconstrained post-order levels; the order is the feasible quantity that moves inventory toward the target. Starting inventory zero and $q_1\ge4\eta$ make all first-period targets in the discrete family feasible.

This is the recursion of \cite[Eq.~(4), p.~358]{Journal} with $T=3$, ACI horizon $n=2$, discount factor $\alpha=1$ and terminal function $f_4=0$. The target and selector agree with its Eq.~(5) and Theorem 2. Starting from future capacities $(\eta,2\eta)$, we compare
\[
 q^{\early}=(0,2\eta),\qquad q^{\late}=(\eta,\eta).
\]
Each removes $\eta$ units from the baseline. For $i\in\{\early,\late\}$, write $J_t^i$ for the scenario objective, $S^i$ for its first-period target set, and $s^i=\min S^i$ for the selected target. To establish that the comparison is admissible under the independent-capacity model of \cite{Journal}, let $Q_2$ have support $\{0,\eta\}$ and $Q_3$ have support $\{\eta,2\eta\}$, all with positive probabilities and independent of demands and each other. All three vectors are possible observed ACI states of this process.

\subsection{Two exact reversals}
Set $p=1/2$, $h=\eta=1$. Substitution in \eqref{eq:bellman} yields the following targets and minimized first-period objectives. The slope calculations establishing uniqueness appear in Appendix~\ref{app:slopes}.
\begin{table}[H]
\centering
\caption{Exact target $s_1$ and minimized cost $J_1(s_1)$ under two equal capacity losses.}\label{tab:examples}
\begin{tabular}{lrrrr}
\toprule
&\multicolumn{2}{c}{$b=1/2$}&\multicolumn{2}{c}{$b=5$}\\
\cmidrule(lr){2-3}\cmidrule(lr){4-5}
Future capacity & Target & Cost & Target & Cost\\\midrule
Baseline $(1,2)$ &0&$15/8$&2&$9/2$\\
Early loss $(0,2)$ &0&$9/4$&2&$6$\\
Late loss $(1,1)$ &1&$33/16$&3&$21/4$\\\bottomrule
\end{tabular}
\end{table}
For $b=1/2$, convexity and the adjacent slopes establish the unique minima at 0 and 1. This violates both temporal assertions in \cite{Journal}. For $b=5$, the target again rises, from 2 to 3, when the loss is moved later. In both cases the earlier loss has the larger optimized cost (Table~\ref{tab:examples}). The examples lie on opposite sides of the one-period critical cost ratio. The following classification determines all cost ratios for which the target order reverses.

\section{Classification of the temporal comparison}\label{sec:main}
\subsection{Target classification}
Write $r=b/h$, $u=1-p$, $c=3-3p+p^2$ and $r_{\mathrm m}=u/p$. The latter is the myopic critical ratio at which the expected one-period cost is flat between 0 and $2\eta$. Define six thresholds, indexed in increasing order within the two cost regimes:
\begin{align}\label{eq:thresholds}
 \tau_1&=\frac{uc}{p(5-3p+p^2)},&
 \tau_2&=\frac{uc}{p(4-2p+p^2)},&
 \tau_3&=\frac{u(3-p)}{3p},\notag\\
 \tau_4&=\frac{3-p-p^2}{2p^2},&
 \tau_5&=\frac{3-2p+p^2-p^3}{p^2(1+p)},&
 \tau_6&=\frac{3-p+p^2-p^3}{p^3}.
\end{align}
The following operator reduces the comparison to the slopes of convex objectives. It will also allow us to establish sensitivity independently of the temporal ranking.

\begin{lemma}[Capacity-constrained minimization]\label{lem:clip}
Let $J$ be a finite coercive convex function on $\mathbb R$ with minimizer interval $[s_-,s_+]$. For $q\ge0$, define $T_qJ(x)=\min_{0\le z\le q}J(x+z)$. Its minimizer interval is $[s_--q,s_+]$. The associated derivative operator $\Phi_q$ is defined by
\begin{equation}\label{eq:clip}
 (T_qJ)'_+(x)=\Phi_q J'_+(x)
 =\max\{J'_+(x),0\}+\min\{J'_+(x+q),0\}.
\end{equation}
\end{lemma}
\begin{proof}
Below $s_--q$, the objective is minimized at $x+q$; above $s_+$ it is minimized at $x$. Between these points the feasible interval intersects $[s_-,s_+]$, so the value is constant. Taking right derivatives on these regions, including their boundary points, gives \eqref{eq:clip}.
\end{proof}

\Needspace{10\baselineskip}
\begin{theorem}\label{thm:phase}
For every $p\in(0,1)$,
\[
0<\tau_1<\tau_2<\tau_3<r_{\mathrm m}<\tau_4<\tau_5<\tau_6.
\]
The complete first-period target sets are
\begin{equation}\label{eq:targetsets}
 \eta^{-1}S^{\early}=\begin{cases}
 \{0\},&r<\tau_2,\\ [0,2],&r=\tau_2,\\
 \{2\},&\tau_2<r<\tau_5,\\{}[2,4],&r=\tau_5,\\\{4\},&r>\tau_5,
 \end{cases}\qquad
 \eta^{-1}S^{\late}=\begin{cases}
 \{0\},&r<\tau_1,\\{}[0,1],&r=\tau_1,\\\{1\},&\tau_1<r<\tau_3,\\
 [1,2],&r=\tau_3,\\\{2\},&\tau_3<r<\tau_4,\\{}[2,3],&r=\tau_4,\\
 \{3\},&\tau_4<r<\tau_6,\\{}[3,4],&r=\tau_6,\\\{4\},&r>\tau_6.
 \end{cases}
\end{equation}
In particular both targets are uniquely $2\eta$ at $r=r_{\mathrm m}$.
\end{theorem}
\begin{proof}
Normalize $h=\eta=1$. The proof has three steps. First, Lemma~\ref{lem:clip} expresses the derivative after capacity-constrained minimization in terms of the original convex objective. Second, applying this operator backward through the two recourse stages gives the slopes on the intervals between demand and capacity breakpoints, separately for $r<r_{\mathrm m}$ and $r>r_{\mathrm m}$. Their zero crossings are precisely the six thresholds in \eqref{eq:thresholds}; the detailed expressions appear in Appendix~\ref{app:slopes}. Third, the explicit positive threshold differences in that appendix establish their ordering, and evaluation at $r=r_{\mathrm m}$ resolves the remaining case. Since the objectives are convex and piecewise linear, the sign changes identify global minimizers, and zero slopes identify entire optimal intervals. Rescaling gives \eqref{eq:targetsets}.
\end{proof}

\begin{figure}[H]
\centering\includegraphics[width=.96\linewidth]{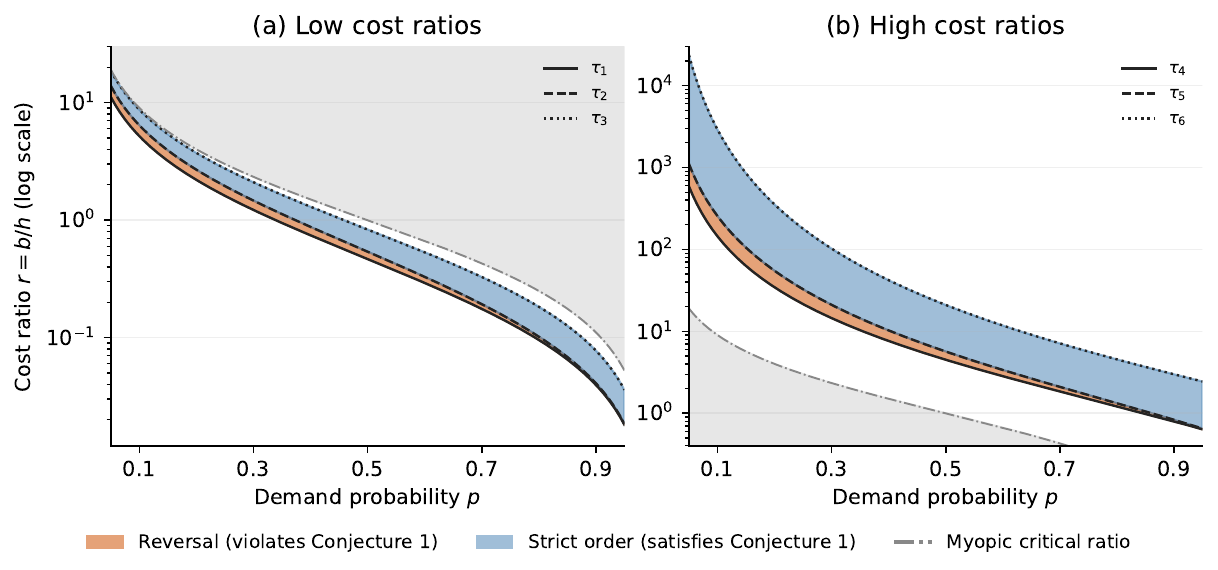}
\caption{Temporal ordering for unique first-period targets. Orange regions violate the conjecture; blue regions satisfy it strictly. White regions within each cost regime have equal targets; gray regions fall outside the regime shown in that panel. The gray dash-dotted line is the myopic critical ratio. Boundary target intervals are given in Theorem~\ref{thm:phase}. Panel (a) displays the low-cost thresholds $\tau_1,\tau_2,\tau_3$; panel (b) displays $\tau_4,\tau_5,\tau_6$. Each uses its own logarithmic vertical range for readability.}\label{fig:phase}
\end{figure}
For unique targets the difference $s^{\late}-s^{\early}$ equals $\eta$ on $(\tau_1,\tau_2)$ and $(\tau_4,\tau_5)$, equals $-\eta$ on $(\tau_2,\tau_3)$ and $(\tau_5,\tau_6)$, and vanishes in the remaining open intervals. Figure~\ref{fig:phase} displays the two regimes. Although the orange bands may look narrow on a logarithmic axis, each is an open nonempty interval for every $p\in(0,1)$. The relative widths in the two regimes are
\[
 \frac{\tau_2}{\tau_1}=1+\frac{u}{4-2p+p^2},\qquad
 \frac{\tau_5}{\tau_4}=1+\frac{uc}{(1+p)(3-p-p^2)}.
\]
Both ratios exceed one for every $p\in(0,1)$ and tend to one as $p$ tends to one.

For the fair-demand case $p=1/2$, the six thresholds are $7/15,7/13,5/6,9/2,17/3,21$, respectively. In particular the conjectured weak order holds outside the two reversal intervals whenever the targets are unique; at the boundaries the full sets in \eqref{eq:targetsets} determine the comparison. Under the smallest-minimizer convention the reversal regions, including their appropriate endpoints, are exactly $(\tau_1,\tau_2]$ and $(\tau_4,\tau_5]$. For a comparison at the first period, three periods are the smallest horizon permitting two distinct future loss locations: with two periods there is only one future capacity to reduce.

\subsection{Cost--target separation}
\begin{proposition}[Optimized cost in the reversal regions]\label{prop:cost}
Let $C^{\early}=\min_y J_1^{\early}(y)$ and $C^{\late}=\min_y J_1^{\late}(y)$. For every $0<p<1$, their difference on the two strict reversal intervals is
\begin{equation}\label{eq:costgap}
\frac{C^{\early}-C^{\late}}{h\eta}=
\begin{cases}
\displaystyle \frac{pu c}{5-3p+p^2}+p(5-2p+p^2)(r-\tau_1),&\tau_1<r<\tau_2,\\[5pt]
\displaystyle \frac{pc}{2}+p^2(2+p)(r-\tau_4),&\tau_4<r<\tau_5.
\end{cases}
\end{equation}
In particular, $C^{\early}>C^{\late}$ throughout both regions, even though $s^{\early}<s^{\late}$.
\end{proposition}
\begin{proof}
Normalize $h=\eta=1$. Theorem~\ref{thm:phase} gives the early and late initial targets $(0,1)$ in the low-cost region and $(2,3)$ in the high-cost region. The terminal target is respectively 0 or 2. In the late-loss scenario the period-2 target may switch between 0 and 1, or between 2 and 3, as established in Appendix~\ref{app:slopes}. At the reachable period-2 inventories (respectively $-1,1$ or $1,3$), either target in the relevant pair induces the same feasible order: one unit at the lower state and zero at the upper state. Hence the intermediate target switch does not change the cost of the chosen first-period decision.
Taking expectations over the eight demand histories using these orders and the terminal target gives
\begin{align*}
\text{Low:}\quad C^{\early}&=2pr(4+p),\\
 C^{\late}&=3-6p+4p^2-p^3+pr(3+4p-p^2),\\[2pt]
\text{High:}\quad C^{\early}&=2(3-4p+p^3)+2p^2(1+p)r,\\
 C^{\late}&=9u+p^3(1+r).
\end{align*}
Subtracting, then substituting $r=\tau_1$ and $r=\tau_4$ into the respective affine differences, yields \eqref{eq:costgap}. Both intercepts and both slopes in that equation are positive. Rescaling restores the factor $h\eta$.
\end{proof}
Even within a region where the conjectured target order holds strictly, the cost order can change. With $p=1/2$ and $h=\eta=1$, both $r=7/10$ and $r=4/5$ lie in $(\tau_2,\tau_3)$ and yield targets $(s^{\early},s^{\late})=(2,1)$. Direct Bellman evaluation gives respective cost differences $7/80$ and $-1/40$. Thus the same strict target order is compatible with either cost order.

\subsection{Economic mechanism} Additional initial stock is valuable when it changes the inventory carried through a binding replenishment constraint. It can also be redundant when later production would restore the same post-order stock. Moving a loss from period 3 to period 2 can raise costs incurred before period 3 while leaving more production available in period 3. The additional terminal capacity can substitute for initial stock. The history $D_1=0,D_2=2$ illustrates this mechanism on both sides of the myopic critical ratio.

In the low-cost example of Table~\ref{tab:examples}, the period-2 and period-3 targets can both be taken as zero. Compare initial orders of zero and one. After $D_1=0$, the starting inventory in period 2 is already at or above the zero target, so neither scenario orders. The extra unit costs one in period 1 and saves $1/2$ in period 2. With the early loss, terminal capacity two clears either backlog, leaving no additional terminal saving. With the late loss, terminal capacity one leaves a backlog of one before demand if initial stock is zero, and no backlog if it is one. The extra unit saves another $1/2$ regardless of final demand. The conditional incremental total cost is thus $1/2$ for the early loss and zero for the late loss.

For $b=5$, the subsequent targets can both be taken as two. Compare initial orders of two and three on the same history. Since $D_1=0$, neither scenario needs to order in period 2 to reach its target of two. The extra unit costs one in each of the first two periods. The early-loss scenario can reach the terminal target two from either inventory state. In the late-loss scenario the two pre-demand terminal stocks are one and two. The additional unit increases terminal holding cost by one if $D_3=0$ but saves backlog cost five if $D_3=2$. Its conditional expected terminal saving is two, so its conditional incremental total cost is zero for the late loss and two for the early loss. In both examples the other three $(D_1,D_2)$ histories have identical incremental costs in the two scenarios. Averaging over the histories gives marginal costs $1/16,-1/16$ in the low-cost example and $1/4,-1/4$ in the high-cost example.

The construction makes this difference possible because a positive demand of $2\eta$ consumes twice the capacity loss, while the remaining total capacity $2\eta$ is either concentrated at the terminal period or split into two smaller replenishments. The timing of the second demand then determines whether the final capacity can undo the shortage. This conditional substitution between current stock and later production explains why the more expensive disruption need not justify greater initial inventory. Proposition~\ref{prop:cost} establishes the cost contrast over both reversal regions, while Theorem~\ref{thm:phase} identifies exactly where the marginal cost changes sign.

\section{Sensitivity and structural stability}\label{sec:sensitivity}
The derivation of the coordinate sensitivity bound in Appendix A of \cite{Journal} uses Conjecture 1 and a first-coordinate lemma. We establish the bound without that premise, then check whether the reversals depend on discrete demand.

\subsection{A sharp coordinate bound}
\begin{proposition}\label{prop:sensitivity}
Consider a finite-horizon zero-order-cost backlogged system with positive linear holding and backlog costs, independent nonnegative demands satisfying $\EE D_t<\infty$, and fixed known future capacities $q$. At every stage use the smallest unconstrained minimizer $s(q)=\min\argmin_y J(y)$ as the target. For any future coordinate $i$ and $0\le\eta\le q_i$, with $e_i$ the $i$th unit vector,
\[
0\le s(q-\eta e_i)-s(q)\le\eta.
\]
Consequently $|s(q-\eta e_i)-s(q-\eta e_j)|\le\eta$ for two equal coordinate losses. For every $p\in(0,1)$ the signed difference $s(q-\eta e_i)-s(q-\eta e_j)$ attains both $\eta$ and $-\eta$ in the family of Theorem~\ref{thm:phase}.
\end{proposition}
\begin{proof}
For a nondecreasing derivative $g$, \eqref{eq:clip} yields
\[
 \Phi_qg(x-\eta)\le\Phi_{q-\eta}g(x)\le\Phi_qg(x).
\]
The derivative operator $\Phi_q$ associated with capacity-constrained minimization preserves derivative order and commutes with translations. Backward expectation preserves the sandwich, and adding the same immediate convex cost also preserves it because the immediate derivative is nondecreasing. Finite-horizon linear costs give bounded slopes, justifying the interchange of right derivative and expectation. Coercivity follows from the positive immediate holding and backlog costs and finite first moments. The smallest minimizer is the first point where the right derivative is nonnegative. Its location therefore moves by an amount in $[0,\eta]$. Both loss targets lie in the same interval $[s(q),s(q)+\eta]$, proving the spread bound. The four strict regions in Theorem~\ref{thm:phase} attain both signs.
\end{proof}
At a flat minimum the common smallest-minimizer convention fixes the target unambiguously. This is essential for the location comparison. The result separates a uniform bound on the magnitude of the response from any ranking by time.

\subsection{Stability under continuous demand perturbations}\label{sec:robust}
The sharpness examples use two-point demand, so their objective functions have exact kinks. We establish local structural stability by perturbing demand with independent continuous noise and bounding the displacement of every optimizer.
\begin{proposition}\label{prop:continuous}
Replace the fair two-point demands in either example of Table~\ref{tab:examples} by iid $D_t^\epsilon=2B_t+U_t$, where $U_t$ is uniform on $[0,\epsilon]$ with the noises mutually independent and independent of the Bernoulli variables. For $0<\epsilon<1/384$ in the low-cost case and $0<\epsilon<1/480$ in the high-cost case, every early-loss minimizer is strictly smaller than every late-loss minimizer. In particular both reversals hold at $\epsilon=10^{-3}$.
\end{proposition}
\begin{proof}
Let $M=\max(h,b)$. The loss is $M$-Lipschitz and a continuation with $H$ stages is at most $HM$-Lipschitz. Coupling the two demand processes and using nonexpansiveness of capacity-constrained minimization gives $\sup_y|J_1^\epsilon(y)-J_1^0(y)|\le6M\epsilon$. If $J_1^0(y)-\min J_1^0\ge\mu|y-y^*|$, each perturbed minimizer is within $12M\epsilon/\mu$ of $y^*$. The adjacent slopes give $\mu\ge1/16$ for both low-cost scenarios, where $M=1$, and $\mu\ge1/4$ for both high-cost scenarios, where $M=5$. Their respective radii are at most $192\epsilon$ and $240\epsilon$. In each example the unperturbed targets differ by one, so the claimed strict bounds make the two neighborhoods disjoint and preserve their order.
\end{proof}

The perturbed distributions are mixtures of two narrow uniform densities. The result establishes local stability around the two explicit demand laws; the ordering under other continuous laws, including unimodal distributions, is a separate question. Having established that the target reversal persists under perturbation, we return to the two-point family to quantify the cost of enforcing the conjectured order.

\section{Decision costs of the temporal ordering}\label{sec:information}
In the two-point-demand family, suppose a single loss of $\eta$ is certain and its location is early with probability $\theta\in(0,1)$ and late otherwise, independently of demand. Consider a planner who observes the loss location before the first order but restricts the two location-dependent targets $y^{\early},y^{\late}$ to satisfy the conjectured ordering $y^{\early}\ge y^{\late}$. The planner otherwise optimizes each first-period order and uses the optimal recourse policy. All optimal targets and their intervening common decisions are feasible by the initial-capacity assumption in Section~\ref{sec:model}. For example, a policy search may enforce the conjecture by restricting its lookup table of location-dependent targets to be monotone in loss timing. The optimization below measures the smallest loss achievable within that restricted policy class.

Let $R^{i}(y)=J_1^i(y)-C^i$ for $i\in\{\early,\late\}$. Define the cost of imposing the order by
\begin{equation}\label{eq:ordercost}
 K(\theta)=\min_{0\le y^{\late}\le y^{\early}\le q_1}\{\theta R^{\early}(y^{\early})+(1-\theta)R^{\late}(y^{\late})\}.
\end{equation}
For comparison, if the location is revealed only after first-period demand and before the period-2 order, the two scenarios must share a common initial target. All subsequent information and physical capacities are the same. The expected value of perfect loss-location information is therefore
\[
 I(\theta)=\min_{0\le y\le q_1}\{\theta R^{\early}(y)+(1-\theta)R^{\late}(y)\}.
\]
Under this two-vector prior, observing either future capacity reveals the location of the single loss. Thus $I$ is the value of revealing the location before the first decision.

\begin{proposition}[Cost of enforcing temporal order]\label{prop:ordercost}
In either strict reversal region of Theorem~\ref{thm:phase}, $K(\theta)=I(\theta)>0$. In all other parameter regions, including the six boundaries, $K(\theta)=0$.
\end{proposition}
\begin{proof}
In a reversal region, project both decisions onto $[s^{\early},s^{\late}]$. This weakly reduces both regrets by convexity and preserves their order. On this interval the early-loss regret is increasing and the late-loss regret is decreasing. For any feasible pair $y^{\early}\ge y^{\late}$, replacing both decisions by any common value between them weakly reduces both regrets. Thus the constrained optimum equals the common-target optimum. Since the two targets are unique and distinct, the common-target optimum has strictly positive expected regret. Outside the strict reversal regions the target sets contain an ordered pair of optimal decisions (at a boundary they intersect), so the constrained optimum is zero.
\end{proof}

The equality has a direct planning interpretation: in the reversal regions, imposing the incorrect order destroys the entire first-period benefit of observing the loss location, even though the information itself is available.

\paragraph{Exact costs.} On each of the four strict-ordering intervals, write $y_- = \min(s^{\early},s^{\late})$ and $y_+=\max(s^{\early},s^{\late})$. The targets differ by $\eta$, and both first-period objectives are affine on $[y_-,y_+]$. For $i\in\{\early,\late\}$ define the normalized slope
\[
 a_i=\frac{J_1^i(y_+)-J_1^i(y_-)}{h\eta}.
\]
The slope is positive for the scenario with the smaller target and negative for the one with the larger target. Convexity restricts a common optimum to this interval, where the prior-weighted regret is affine. Evaluating at the two endpoints yields
\begin{equation}\label{eq:evpi}
 I(\theta)=h\eta\min\{\theta|a_{\early}|,(1-\theta)|a_{\late}|\}.
\end{equation}
For $(\tau_1,\tau_2)$ and $(\tau_2,\tau_3)$ the early slope is the first line of \eqref{eq:elow}; the late slope is the derivative on $(0,1)$ and $(1,2)$, respectively, in \eqref{eq:dlow}. For $(\tau_4,\tau_5)$ and $(\tau_5,\tau_6)$ the early slope is the derivative on $(2,4)$ in \eqref{eq:ehigh}; the late slope is the derivative on $(2,3)$ and $(3,4)$, respectively, in \eqref{eq:dhigh}. These formulas in Appendix~\ref{app:slopes} use normalized units. Outside these four intervals the target sets intersect, so $I(\theta)=0$. By Proposition~\ref{prop:ordercost}, the same formula gives $K(\theta)$ in the two reversal regions, whereas $K(\theta)=0$ in the two strictly consistent regions.

For a fair prior and $p=1/2$, $h=\eta=1$, the first two rows of Table~\ref{tab:policycost} revisit the reversals in Table~\ref{tab:examples}. The third uses $r=7/10$, which lies between $\tau_2=7/13$ and $\tau_3=5/6$, and satisfies the conjecture. In the two reversal examples the order restriction incurs the full information value; in the consistent example it costs nothing even though withholding information is costly. The percentages use the expected optimized cost with full information as the common denominator for each example.
\begin{table}[htbp]
\centering
\caption{Fair-prior costs ($p=1/2$, $h=\eta=1$). Percentages are relative to the expected full-information cost.}\label{tab:policycost}
\begin{tabular}{llrrr}\toprule
Regime & $r$ & Full-info cost & $I(1/2)$ & $K(1/2)$\\\midrule
Reversal (low) & $1/2$ & $69/32$ & $1/32$ (1.45\%) & $1/32$ (1.45\%)\\
Reversal (high) & $5$ & $45/8$ & $1/8$ (2.22\%) & $1/8$ (2.22\%)\\
Consistent (low) & $7/10$ & $413/160$ & $1/10$ (3.87\%) & $0$ (0\%)\\\bottomrule
\end{tabular}
\end{table}

\section{Conclusion}\label{sec:conclusion}
For the three-period two-point-demand family, six explicit thresholds determine both optimal inventory target sets for every positive cost ratio and nondegenerate demand probability. In each of the two strict reversal regions the earlier capacity loss has strictly higher optimized cost despite its lower initial target. The coordinate sensitivity bound nevertheless holds independently of the temporal comparison and is sharp in both directions. Small continuous demand perturbations preserve the two explicit reversals.

In the fixed-loss experiment, enforcing the conjectured ordering on otherwise informed decisions eliminates the entire first-period value of the loss-location information in the two strict reversal regions. Outside those regions the constraint has no cost. The distinction between the severity of a disruption and the optimal response to it is therefore consequential for decision rules, not merely for the ranking of two stock levels. Extending the cost and target classifications to longer horizons and general demand distributions, and evaluating partial forecasts under independently distributed future capacities, remain natural directions.

\paragraph{Correspondence.} Qichao Wang, \href{mailto:cuilijun42@gmail.com}{cuilijun42@gmail.com}.

\paragraph{Code and data availability.} Reproduction scripts are available at \url{https://github.com/dongxuelian2/computational-economic-planning}, under \path{scripts/}: \path{reproduce_aci_counterexample.py}, \path{aci_piecewise.py}, \path{audit_aci_full_phase.py}, and \path{check_aci_scenario_lp.py}. The first two implement the exact counterexamples and piecewise-linear recursion; the latter scripts check the classification and the scenario-tree optimization. The script \path{check_aci_cost_order.py} independently evaluates the cost differences and both decision problems using exact rational Bellman recursion and direct enumeration of feasible pairs at integer knots. No external empirical data are required.

\appendix
\Needspace{6\baselineskip}
\section{Bellman slope calculations}\label{app:slopes}
Normalize $h=\eta=1$. Let $G(y)=\EE\ell(y-D)$ and $m=u-pr$. Its slopes on intervals separated by $0,2$ are $-r,m,1$. Every objective in \eqref{eq:bellman} in this two-point family is finite, convex and piecewise linear. We use Lemma~\ref{lem:clip} at each stage. The intermediate recourse ratios are $\rho_-=u(2-p)/[p(3-p)]$ and $\rho_+=2/p^2-1$. All derivatives below are in normalized units and hence equal the original slopes divided by $h$.
\subsection{Low cost ratios}
Here $r<r_{\mathrm m}$, $m>0$ and $G$ has unique target zero. With $q_3=2$, $V_3$ has slopes $-r,0,m,1$ at breakpoints $-2,0,2$. The early-loss period-2 objective has slopes
\[
 -(1+p)r\quad(-2,0),\qquad m(2-p)\quad(0,2),\qquad 2-p+pm\quad(2,4).
\]
As $q_2=0$, $V_2=J_2$, and
\begin{align}\label{eq:elow}
 (J_1^{\early})'|_{(0,2)}&=uc-pr(4-2p+p^2)=p(4-2p+p^2)(\tau_2-r),\notag\\
 (J_1^{\early})'|_{(2,4)}&=c+pm(3-2p)>0.
\end{align}
The objective decreases strictly below zero; hence \eqref{eq:elow} proves the early-loss alternatives at $\tau_2$.

For the late loss, $q_3=1$ and terminal continuation slopes $-r,0,m,1$ have breakpoints $-1,0,2$. Let
\[
 v_0=u(2-p)-pr(3-p).
\]
Period-2 slopes on $(0,1),(1,2),(2,3),(3,4)$ are $v_0,m(2-p),2-p+pm,2-p+pm$. The target is zero below $\rho_-$, one above $\rho_-$, and $[0,1]$ at equality. Capacity-constrained minimization with $q_2=1$ gives
\begin{align}\label{eq:dlow}
 (J_1^{\late})'|_{(0,1)}&=m+u\max(v_0,0)-p(1+p)r,\notag\\
 (J_1^{\late})'|_{(1,2)}&=mc+p\min(v_0,0).
\end{align}
Below $\rho_-$ these are $p(5-3p+p^2)(\tau_1-r)$ and $mc>0$; above $\rho_-$ they are $u-pr(2+p)$ and $3p(\tau_3-r)$. The identities
\[
 \rho_--\tau_1=\frac{u(1+p-p^2)}{p(3-p)(5-3p+p^2)}>0,\qquad
 \tau_3-\rho_-=\frac{uc}{3p(3-p)}>0
\]
show that the two current slopes cross zero only at $\tau_1$ and $\tau_3$. At $\rho_-$, the slope on $(0,1)$ is negative and that on $(1,2)$ is positive; both branch formulas agree. The derivative below zero is negative and that on $(2,3)$ is positive. Convexity gives the low-cost late-loss classification.

\subsection{High cost ratios}
Now $r>r_{\mathrm m}$, $m<0$ and $G$ has unique target 2. For the early loss, $V_3$ has slopes $-r,m,0,1$ at breakpoints $-2,0,2$. The period-2 target is 2, and the current slopes are
\begin{equation}\label{eq:ehigh}
 (J_1^{\early})'(y)=\begin{cases}
 m(2+p-2p^2)-p(1+p)r<0,&0<y<2,\\
 p^2(1+p)(\tau_5-r),&2<y<4,\\3-p>0,&4<y<6.
 \end{cases}
\end{equation}
The middle slope gives the early-loss transition at $\tau_5$, including $[2,4]$ at equality.

For the late loss, $V_3$ has slopes $-r,m,0,1$ at $-1,1,2$. Period-2 slopes on the four consecutive unit intervals from 0 to 4 are
\[
 v_0,\quad m(1+p),\quad v_2=2-p^2(1+r),\quad 2-p.
\]
Its target changes from 2 to 3 at $\rho_+$; at equality it is $[2,3]$. The current middle slopes are
\begin{equation}\label{eq:dhigh}
 \begin{array}{c|cc}
 & (J_1^{\late})'\text{ on }(2,3)&(J_1^{\late})'\text{ on }(3,4)\\\hline
 r\le\rho_+&2p^2(\tau_4-r)&c\\
 r\ge\rho_+&1+p-p^3-p^2(1+p)r&p^3(\tau_6-r).
 \end{array}
\end{equation}
The formulas agree at $\rho_+$, where the $(2,3)$ slope is $-1-p+p^2<0$. The $(1,2)$ slope is negative: below $\rho_+$ it is $m+pv_0<0$, and above $\rho_+$ there is the additional nonpositive term $uv_2$. The $(4,5)$ slope is positive in both regimes. Since $\tau_4<\rho_+<\tau_6$, the current target crosses at $\tau_4$ and $\tau_6$.

\subsection{Threshold order and the myopic boundary}
The following positive differences establish the remaining comparisons:
\begin{align*}
 \tau_2-\tau_1&=\frac{u^2c}{p(5-3p+p^2)(4-2p+p^2)},&
 \tau_3-\tau_2&=\frac{u(3-p+2p^2-p^3)}{3p(4-2p+p^2)},\\
 r_{\mathrm m}-\tau_3&=u/3,&\tau_4-r_{\mathrm m}&=\frac{c}{2p^2},\\
 \tau_5-\tau_4&=\frac{uc}{2p^2(1+p)},&
 \tau_6-\tau_5&=\frac{3-p+2p^2-p^3}{p^3(1+p)}.
\end{align*}
Also $\rho_+-\tau_4=(1+p-p^2)/(2p^2)>0$ and $\tau_6-\rho_+=c/p^3>0$. At $r=r_{\mathrm m}$, terminal and recourse optima may be intervals, but the early current slopes adjacent to 2 are $-(1+p)u<0$ and $c>0$; the late adjacent slopes are $-pu<0$ and $c>0$. Both current targets are uniquely 2.

\end{document}